\documentclass[reqno]{amsart}
\usepackage{amsmath,amssymb,amsfonts,amsthm,longtable}

\newtheorem{theorem}{Theorem}
\newtheorem{lemma}[theorem]{Lemma}
\newtheorem{proposition}[theorem]{Proposition}
\newtheorem{corollary}[theorem]{Corollary}
\newtheorem{conjecture}[theorem]{Conjecture}
\newtheorem{question}[theorem]{Question}
\theoremstyle{definition}

\newtheorem{remark}[theorem]{Remark}
\newtheorem{example}{Example}

\newcommand{\bij}{\Gamma}
\newcommand{\biject}{\bij}

\newcommand{\invbij}{\bij^{-1}}
\newcommand{\asc}{\mathsf{asc}}

\newcommand{\last}{\mathsf{last}}
\newcommand{\comp}{\mathsf{comp}}
\newcommand{\rmax}{\mathsf{rmax}}
\newcommand{\blocks}{\mathsf{blocks}}
\newcommand{\zeros}{\mathsf{zeros}}

\newcommand{\Asc}[1]{{\mathcal{A}}_{#1}}

\newcommand{\tpt}{$(\mathbf{2+2})$}

\renewcommand{\matrix}{\left(\begin{array}{cc}}
\newcommand{\twomatrix}{\left(\begin{array}{cc}}
\newcommand{\threematrix}{\left(\begin{array}{ccc}}
\newcommand{\fourmatrix}{\footnotesize\left(\begin{array}{cccc}}
\newcommand{\fivematrix}{\footnotesize\left(\begin{array}{ccccc}}
\newcommand{\sixmatrix}{\footnotesize\left(\begin{array}{cccccc}}
\newcommand{\sevenmatrix}{\footnotesize\left(\begin{array}{ccccccc}}
\newcommand{\ematrix}{\end{array}\right)}
\newcommand{\ourclass}[1]{\mathsf{Int}_{#1}}

\newcommand{\bidiagclass}[1]{\mathsf{Bi}_{#1}}
\newcommand{\ourindex}[1]{\mathsf{index}({#1})}
\newcommand{\ourvalue}[1]{\mathsf{value}({#1})}
\newcommand{\ourdim}[1]{\mathsf{dim}({#1})}
\newcommand{\rowsum}[2]{\mathsf{rowsum}_{#1}(#2)}
\newcommand{\colsum}[2]{\mathsf{colsum}_{#1}(#2)}
\newcommand{\ascset}{\mathfrak{asc}}
\newcommand{\addi}{\textsf{Add1}}
\newcommand{\addii}{\textsf{Add2}}
\newcommand{\addiii}{\textsf{Add3}}
\newcommand{\flip}{\mathsf{flip}}
\newcommand{\ds}{\displaystyle}

\newcommand{\indicator}[1]{\mathbf{1}(#1)}

\newcommand{\subi}{\textsf{Rem1}}
\newcommand{\subii}{\textsf{Rem2}}
\newcommand{\subiii}{\textsf{Rem3}}

\newcommand{\addto}[1]{\phi(#1)}

\title[Ascent sequences and upper triangular matrices]{Ascent sequences and upper triangular matrices containing non-negative integers}
\author[M. Dukes]{Mark Dukes}
\author[R. Parviainen]{Robert Parviainen}
\thanks{Both authors were supported by grant no. 090038011 from the Icelandic Research Fund.}
\address{M. Dukes: Science Institute, University of Iceland, 107 Reykjav\'ik, Iceland}
\address{R. Parviainen: The Mathematics Institute, School of Computer Science, Reykjav\'ik University, 103 Reykjav\'ik, Iceland}

\begin{document}
\maketitle
\begin{abstract}
This paper presents a bijection between ascent sequences and upper triangular matrices whose non-negative entries are such that
all rows and columns contain at least one non-zero entry.
We show the equivalence of several natural statistics on these structures under this bijection and 
prove that some of these statistics are equidistributed. 
Several special classes of matrices are shown to have simple formulations in terms of ascent sequences.
Binary matrices are shown to correspond to ascent sequences with no two adjacent entries the same.
Bidiagonal matrices are shown to be related to order-consecutive set partitions and a simple condition on the ascent sequences
generate this class.
\end{abstract}

\section{Introduction}
Let $\ourclass{n}$ be the collection of upper triangular matrices 
with non-negative integer entries which sum to $n \in \mathbb{N}$ such that 
all rows and columns contain at least one non-zero entry.
For example,
$$\ourclass{3} = \left\{(3) , 
	\twomatrix 1 & 1 \\ 0 & 1 \ematrix ,
	\twomatrix 2 & 0 \\ 0 & 1 \ematrix ,
	\twomatrix 1 & 0 \\ 0 & 2 \ematrix ,
	\threematrix 1 & 0 & 0 \\ 0 & 1 & 0 \\ 0 & 0 & 1 \ematrix
	\right\}.
$$

We use the standard notation $[a,b]$ for the interval of integers $\{a, a+1, \dots, b\}$ and define $[n]=[1,n]$.
Given a sequence of integers $y=(y_1,\ldots , y_n)$, we say that $y$ has an {\it{ascent}} at position $i$ if $y_i<y_{i+1}$.
The number of ascents of $y$ is denoted by $\asc(y)$.
Let $\Asc{n}$ be the collection of {\it{ascent sequences of length $n$}}:
$$\Asc{n} = \{(x_1,\ldots , x_n): 
	x_i \in[0,1+\asc(x_1,\ldots , x_{i-1})] , 
	\mbox{ for all }1<i\leq n \},$$
where $x_1:=0$ and $\asc(x_1):=0$.
For example,
$$\Asc{3}=\{
	(0,0,0), \, (0,0,1), \, (0,1,0), \, (0,1,1), (0,1,2)
	\}.$$

These sequences were introduced in the recent paper by Bousquet-M\'elou et al. ~\cite{bcdk}
and were shown to unify three combinatorial structures: \tpt -free posets, a class of pattern avoiding permutations
and a class of involutions that are sometimes termed {\it{chord diagrams}}.
This paper complements the results of \cite{bcdk} by presenting a fourth structure, the matrices in $\ourclass{n}$, 
that can be encoded by an ascent sequence of length $n$.
To this end we have attempted to use notation that is indicative of the transformations and operations in the original paper \cite{bcdk}.

The class of matrices we study here have been touched upon in the literature before. 
The binary case is known to encode a subclass of interval orders 
(the full class of interval orders are in bijection with \tpt -free posets), see Fishburn \cite{F85}. 
Mitas \cite{M91} used our class of matrices to study the jump number problem on interval orders, 
but without a formal statement or proof of any bijection, and without studying further properties of the relation.

In section 2 we present a bijection $\biject$ from 
matrices in $\ourclass{n}$ to ascent sequences in $\Asc{n}$.
In section 3 we show how statistics on both of these structures are related under $\biject$
and prove that some of the statistics are equidistributed.
Section 4 looks at properties of restricted sets of matrices and ascent sequences which give rise
to interesting structures, order-consecutive set partitions being one example.
We end with some open problems in section \ref{sec:Q}.

\section{Upper triangular matrices}
In this section we will define a removal and an addition operation on matrices in $\ourclass{n}$ that are essential for the bijection.
These operations have the effect of decreasing (resp. increasing) the sum of the entries in a matrix by 1.

Given $A \in \ourclass{n}$ let $\ourdim{A}$ be the number of rows in the matrix $A$.
Furthermore, let $\ourindex{A}$ be the smallest value of $i$ such that $A_{i,\ourdim{A}}>0$ and define
$\ourvalue{A} := A_{\ourindex{A},\ourdim{A}}$.

Consider the following operation $f$ on a given matrix $A \in \ourclass{n}$.
\begin{enumerate}
\item[(\subi)] If $\ourvalue{A}>1$, or if $\ourvalue{A}=1$ and $\ourindex{A}<\ourdim{A}$, and there is at least one other positive entry in row $\ourindex{A}$, then let $f(A)$ be the matrix $A$ with the entry $A_{\ourvalue{A}, \ourdim{A}}$ reduced by 1. 
\item[(\subii)] If $\ourvalue{A}=1$ and $\ourindex{A}=\ourdim{A}$ let $f(A)$ be the matrix $A$ with row $\ourdim{A}$ and column $\ourdim{A}$ removed.
\item[(\subiii)]
If $\ourvalue{A}=1$, $\ourindex{A}<\ourdim{A}$, and all other entries in row $\ourindex{A}$ are 0, then we form $f(A)$ in the following way.
Let $A_{i,\ourdim{A}}=A_{i,\ourindex{A}}$ for all $1\leq i \leq \ourindex{A}-1$.
Now simultaneously delete row $\ourindex{A}$ and column $\ourindex{A}$.
Let the resulting $(\ourdim{A}-1)\times (\ourdim{A}-1)$ matrix be $f(A)$.
\end{enumerate}
\begin{example}
Consider the following three matrices:
$$A=\fourmatrix
	1&0&1&0\\
	0&2&0&3\\
	0&0&1&4\\
	0&0&0&2
	\ematrix ;\quad
B=\fourmatrix
	5&1&3&0\\
	0&1&0&0\\
	0&0&1&0\\
	0&0&0&1
	\ematrix ;\quad 
C=\sevenmatrix
	1&0&0&1&0&0&0\\
	0&1&0&1&1&0&0\\
	0&0&1&2&1&1&0\\
	0&0&0&0&0&0&1\\
	0&0&0&0&0&1&0\\
	0&0&0&0&0&0&1\\
	0&0&0&0&0&0&1
	\ematrix .$$
For matrix $A$, rule $\subi$ applies since $\ourvalue{A}=3$ and 
$$f(A) = \fourmatrix
	1&0&1&0\\
	0&2&0&2\\
	0&0&1&4\\
	0&0&0&2
	\ematrix  .$$
For matrix $B$, since $\ourvalue{B}=1$ and $\ourindex{B}=\ourdim{B}=4$ rule $\subii$ applies and 
$$f(B)=\threematrix
	5&1&3\\
	0&1&0\\
	0&0&1
	\ematrix .$$
For matrix $C$, since $\ourvalue{C}=1$, $4=\ourindex{C}<\ourdim{C}=7$, and all other entries in row $\ourindex{C}=4$ are zero, 
then we form $f(C)$ in the following way:
first copy the $\ourindex{C}-1=3$ highest entries in column $\ourindex{C}$ to the top $\ourindex{C}-1=3$ entries in column $\ourdim{C}=7$. These are illustrated in bold in the following matrix:
$$\sevenmatrix
	1&0&0&1&0&0&\bf{1}\\
	0&1&0&1&1&0&\bf{1}\\
	0&0&1&2&1&1&\bf{2}\\
	0&0&0&0&0&0&1\\
	0&0&0&0&0&1&0\\
	0&0&0&0&0&0&1\\
	0&0&0&0&0&0&1
	\ematrix .$$
Next we simultaneously remove column $\ourindex{C}=4$ and row $\ourindex{C}=4$ to get $f(C)$:
$$	\sevenmatrix
	1&0&0&\phantom{2}&0&0&\bf{1}\\
	0&1&0& &1&0&\bf{1}\\
	0&0&1& &1&1&\bf{2}\\
	&&&&&&\\
	0&0&0& &0&1&0\\
	0&0&0& &0&0&1\\
	0&0&0& &0&0&1
	\ematrix 
	\quad
	\Longrightarrow 
	\quad
	f(C)=
	\sixmatrix
	1&0&0&0&0&\bf{1}\\
	0&1&0&1&0&\bf{1}\\
	0&0&1&1&1&\bf{2}\\
	0&0&0&0&1&0\\
	0&0&0&0&0&1\\
	0&0&0&0&0&1
	\ematrix .
	$$
\end{example}

\newcommand{\removal}[1]{\psi(#1)}
We now show that the above removal operation yields an upper triangular matrix in $\ourclass{n-1}$.
If $\ourindex{A}=i+1$ and the above removal operation, applied to $A$, gives $f(A)$, then we define $\removal{A}=(f(A),i)$.
Notice that $1\leq\ourindex{A}\leq \ourdim{A}$.

\begin{lemma}
If $n\geq 2$, $A \in \ourclass{n}$ and $\removal{A}=(B,i)$, then $B \in \ourclass{n-1}$. 
\end{lemma}
\begin{proof}
Consider first the case that $\ourvalue{A}>1$. 
Rule \subi\ applies and $B$ is $A$ except that 
$B_{\ourindex{A},\ourdim{A}}=A_{\ourindex{A},\ourdim{A}}-1>0$. 
Thus $B\in \ourclass{n-1}$.

When $\ourvalue{A}=1$ there are three subcases to consider.

First assume that $\ourindex{A}<\ourdim{A}$, 
and that there is a least one other positive entry in row $\ourindex{A}$. 
Rule \subi\ applies and $B$ is $A$ with the 1 at position $(\ourindex{A},\ourdim{A})$ reduced to a 0. 
Since there is another positive entry in the same row and $B_{\ourdim{A},\ourdim{A}}>0$, we have $B \in \ourclass{n-1}$.

Next, if $\ourindex{A}=\ourdim{A}$ then rule \subii\ applies, 
and $B$ is $A$ with row $\ourdim{A}$ and column $\ourdim{A}$ removed.
Since the row $\ourdim{A}-1$ in $A$ has a 0 at position $\ourdim{A}$, 
the value $A_{\ourdim{A}-1,\ourdim{A}-1}>0$. 
Therefore $B \in \ourclass{n-1}$.

Finally, if $\ourindex{A}<\ourdim{A}$, and all other entries in row $\ourindex{A}$ are 0, then rule \subiii\ applies. 
To create $B$, row and column $\ourindex{A}$ are removed from $A$, 
and all positive entries from column $\ourindex{A}$ in $A$ are copied to the last column of $B$. 
Thus no row of zeros or column of zeros is created, and $B \in \ourclass{n-1}$.
\end{proof}


We now define the complementary addition rules for each of the removal steps. Their consistency will be shown later.
Given $A \in \ourclass{n}$ and $m \in [0,\ourdim{A}]$ we construct the matrix $\addto{A,m}$ in the following manner.
\begin{enumerate}
\item[(\addi)] If $0\leq m \leq \ourindex{A}-1$ then let $\addto{A,m}$ be the matrix $A$ with the entry at position $(m+1,\ourdim{A})$ increased by 1.
\item[(\addii)] If $m=\ourdim{A}$ then let $\addto{A,m}$ be the matrix $\left(\begin{array}{cc}A&0\\ 0 & 1 \end{array}\right)$.
\item[(\addiii)]
If $\ourindex{A}\leq m <\ourdim{A} $ then form $\addto{A,m}$ in the following way:

In $A$, insert a new (empty) row between rows $m$ and $m+1$, and insert a new (empty) column between columns $m$ and $m+1$.
Let the new row be filled with all zeros except for the rightmost entry which is 1.
Move each of the entries above this new rightmost one to the new column between columns $m$ and $m+1$ and replace them with zeros.
Finally let all other entries in the new column be zero. The resulting matrix is $\addto{A,m}$.
\end{enumerate}

\begin{example}
Consider the following three matrices:
$$A=\fourmatrix
	1&0&1&0\\
	0&2&0&0\\
	0&0&1&5\\
	0&0&0&1
	\ematrix ; \quad
B=\fourmatrix
	1&5&0&4\\
	0&1&0&3\\
	0&0&1&2\\
	0&0&0&3
	\ematrix ;\quad
C=\sixmatrix
	1&0&0&0&6&0\\
	0&1&0&1&0&7\\
	0&0&1&1&1&2\\
	0&0&0&0&3&0\\
	0&0&0&0&0&1\\
	0&0&0&0&0&1
	\ematrix .$$
In order to form $\addto{A,1}$, since $m=1\leq \ourindex{A}-1 =2$ we see that rule \addi\ applies and 
$$\addto{A,1}=\fourmatrix
        1&0&1&0\\
        0&2&0&1\\
        0&0&1&5\\
        0&0&0&1
        \ematrix .$$
In order to form $\addto{B,4}$, since $m=4=\ourdim{B}$ we see that rule \addii\ applies and 
$$\addto{B,4}=\fivematrix
	1&5&0&4&0\\
	0&1&0&3&0\\
	0&0&1&2&0\\
	0&0&0&3&0\\
	0&0&0&0&1
	\ematrix . $$
In order to form $\addto{C,3}$, since $\ourindex{C}=2 \leq 3<5=\ourdim{C}$ we see that rule \addiii\ applies and we do as follows.
Insert a new empty row and column between rows 3 and 4 and columns 3 and 4 of $C$:
$$
\sevenmatrix
	1&0&0&\phantom{0}&0&6&0\\
	0&1&0& &1&0&7\\
	0&0&1& &1&1&2\\ \\
	0&0&0& &0&3&0\\
	0&0&0& &0&0&1\\
	0&0&0& &0&0&1
	\ematrix .$$
Fill the empty row with all zeros and a rightmost 1, this is highlighted in bold.
Next move the entries above the new {\bf{1}} to the new column and replace them with zeros.
$$
\sevenmatrix
	1&0&0&\phantom{0}&0&6&0\\
	0&1&0& &1&0&7\\
	0&0&1& &1&1&2\\ 
	\bf{0}&\bf{0}&\bf{0}&\bf{0}&\bf{0}&\bf{0}&\bf{1}\\
	0&0&0& &0&3&0\\
	0&0&0& &0&0&1\\
	0&0&0& &0&0&1
\ematrix 
\quad \rightarrow \quad
\sevenmatrix
	1&0&0&{\bf{0}}&0&6&{\bf{0}}\\
	0&1&0&{\bf{7}}&1&0&{\bf{0}}\\
	0&0&1&{\bf{2}}&1&1&{\bf{0}}\\ 
	\bf{0}&\bf{0}&\bf{0}&\bf{0}&\bf{0}&\bf{0}&\bf{1}\\
	0&0&0& &0&3&0\\
	0&0&0& &0&0&1\\
	0&0&0& &0&0&1
	\ematrix .$$
Finally fill the remaining empty positions with zeros to yield $\addto{C,3}$:
$$\addto{C,3} = \sevenmatrix
	1&0&0&{\bf{0}}&0&6&{\bf{0}}\\
	0&1&0&{\bf{7}}&1&0&{\bf{0}}\\
	0&0&1&{\bf{2}}&1&1&{\bf{0}}\\ 
	\bf{0}&\bf{0}&\bf{0}&\bf{0}&\bf{0}&\bf{0}&\bf{1}\\
	0&0&0&{\bf{0}}&0&3&0\\
	0&0&0&{\bf{0}} &0&0&1\\
	0&0&0&{\bf{0}} &0&0&1
	\ematrix .$$
\end{example}

We now show that this addition operation yields another upper triangular matrix where every row and column contain at
least one non-zero entry.

\begin{lemma}
\label{lemma:three}
If $n\geq 2$, $B \in \ourclass{n-1}$, $0\leq i \leq \ourdim{B}$ and $A=\addto{B,i}$, 
then $A \in \ourclass{n}$ and $\ourindex{A} = i+1$. 
\end{lemma}

\begin{proof}
In each of the operations, \addi , \addii\ and \addiii , the sum of the entries of the matrix is increased by exactly 1.
It is straightforward to check that each row and column contains at least one non-zero entry.
The property of being upper-triangular is also preserved. Thus it is clear that $A=\addto{B,i} \in \ourclass{n}$.

It is similarly straightforward to check that $\ourindex{A}=i+1$ in each of the three cases.
\end{proof}

\begin{lemma}
For any $B \in \ourclass{n}$ and integer $i$ such that $0\leq i \leq\ourdim{B}$, we have
$\removal{\addto{B,i}}=(B,i)$. If $n>1$ then we also have $\addto{\removal{B}}=B$.
\end{lemma}

\begin{proof}
First let us denote $A=\addto{B,i}$. 
From Lemma \ref{lemma:three} above $\ourindex{A}=i+1$ and so the removal operation when applied to $A$ 
will yield $\removal{A}=(C,i)$ for some matrix $C$. Thus we need only show that $B=C$ for each of the three cases.

Let us assume that $0\leq i \leq \ourindex{B}-1$. Then $A$ is simply a copy of $B$ with the entry at position
$(i+1,\ourdim{B})$ increased by one. 
Similarly, rule \subi\ applies for $A$ and so $C$ will be the same as $A$ except that the entry at position
$(\ourindex{A},\dim{B})=(i+1,\dim{B})$ is decreased by one. Thus $B=C$.

Assume next that $i=\ourdim{B}$, so that rule \addii\ applies and $A=\twomatrix B&0 \\ 0&1 \ematrix$. 
Since $\ourindex{A}=\ourdim{A}$, rule \subii\ applies and we remove both column and row $\dim{A}$ of $A$ to get $C=(B)$.

If $\ourindex{B}\leq i <\ourdim{B}$ then rule \addiii\ applies.
For this, $B$ must have the following form
$$
B=\left(
	\begin{array}{ccc}
	\mbox{\Huge{$X$}} & \mbox{\Huge{$Y$}} & \begin{array}{c}e_1 \\ \vdots \\ e_i \end{array} \\
	\mbox{\Huge{$0$}} & \mbox{\Huge{$Z$}} & \begin{array}{c}e_{i+1} \\ \vdots \\ e_n \end{array} 
	\end{array}
\right)
$$
where at least one of $\{e_1,\ldots , e_i\}$ is non-zero.
From this we find that
$$
A=\left(
	\begin{array}{ccccc}
	\mbox{\Huge{$X$}} & \begin{array}{c}e_1 \\ \vdots \\ e_i \end{array}  & \mbox{\Huge{$Y$}} & \begin{array}{c}0 \\ \vdots \\ 0 \end{array} \\
	0\cdots 0 &0&0\cdots 0 & 1 \\
	\mbox{\Huge{$0$}} & \begin{array}{c}0 \\ \vdots \\ 0 \end{array}& \mbox{\Huge{$Z$}} & \begin{array}{c}e_{i+1} \\ \vdots \\ e_n \end{array} 
	\end{array}
\right).
$$
Since $\ourindex{A}=i+1$, $\ourvalue{A}=1$ and all other entries in this row are zero, 
the removal operation to be applied is \subiii\ and we find that
$$
C=\left(
	\begin{array}{ccccc}
	\mbox{\Huge{$X$}} & \phantom{\begin{array}{c}e_1 \\ \vdots \\ e_i \end{array}}  & \mbox{\Huge{$Y$}} & \begin{array}{c}e_1 \\ \vdots \\ e_i \end{array} \\
	\phantom{0\cdots 0} &\phantom{0}&\phantom{0\cdots 0} & \phantom{1} \\
	\mbox{\Huge{$0$}} & \phantom{\begin{array}{c}0 \\ \vdots \\ 0 \end{array}}& \mbox{\Huge{$Z$}} & \begin{array}{c}e_{i+1} \\ \vdots \\ e_n \end{array} 
	\end{array}
\right) = B.
$$
The second statement follows by applying a similar analysis of the addition and removal operations.
\end{proof}

We now define a map $\biject$ from $\ourclass{n}$ to $\Asc{n}$ recursively as follows.
For $n=1$ we let $\biject((1))=(0)$.
Now let $n\geq 2$ and suppose that the removal operation, when applied to $A \in \ourclass{n}$, gives $\removal{A}=(B,i)$.
Then the sequence associated with $A$ is $\biject(A) := (x_1,\ldots , x_{n-1},i)$, where $(x_1,\ldots , x_{n-1})=\biject(B)$.

\begin{theorem}\label{thm:biject}
The map $\biject : \ourclass{n} \mapsto \Asc{n}$ is a bijection.
\end{theorem}

\begin{proof}
Since the sequence $\biject (A)$ encodes the construction of the matrix $A$, the map $\biject$ is injective.
We want to prove that the image of $\ourclass{n}$ is the set $\Asc{n}$.
The recursive description of the map $\biject$ tells us that $x=(x_1,\ldots , x_n) \in \biject(\ourclass{n})$ if and only if
\begin{eqnarray}
\label{eq:four}
x'=(x_1,\ldots , x_{n-1}) \in \biject(\ourclass{n-1}) & \mbox{ and } & 0\leq x_n \leq \ourdim{\biject^{-1}(x')}.
\end{eqnarray}
We will prove by induction on $n$ that for all $A \in \ourclass{n}$, with associated sequence $\biject(A)=x=(x_1,\ldots , x_n)$, 
one has 
\begin{eqnarray}
\label{eq:five}
\ourdim{A}=\asc(x) &\mbox{ and }& \ourindex{A}=x_n+1.
\end{eqnarray}
Clearly, this will convert the above description (\ref{eq:four}) of $\biject(A)$ into the definition of ascent sequences,
thus concluding the proof.

So let us focus on the properties (\ref{eq:five}).
They hold for $n=1$. Assume they hold for some $n-1$ with $n\geq 2$,
and let $A=\addto{B,i}$ for $B \in \ourclass{n-1}$.
If $\biject(B)=x' = (x_1,\ldots , x_{n-1})$ then $\biject(A) = x = (x_1,\ldots , x_{n-1},i)$.

Lemma \ref{lemma:three} gives $\ourindex{A}=i+1$ and it follows that
$$
\ourdim{A} = \left\{
	\begin{array}{ll}
	\ourdim{B}=\asc (x')=\asc(x) & \mbox{ if } i \leq x_{n-1},\\
	\ourdim{B}+1=\asc (x')+1=\asc(x) & \mbox{ if } i > x_{n-1}.
	\end{array}
	\right.
$$
The result follows.
\end{proof}

Consider an ascent sequence $x=(x_1, x_2, \ldots, x_n)$. 
The sequence encodes the corresponding matrix $A=\invbij{(x)}$, 
and can be viewed as instructions on how to build up $A$ step by step starting from the single matrix $(1) \in \ourclass{1}$. 
We will elaborate here on this encoding, as it will be used repeatedly. 

Assume that $\bij(A^{(i)})=(x_1, x_2, \ldots, x_i)$. 
For the next step in building up $A=A^{(n)}$ according to $x$, there are three cases depending on the value $x_{i+1}$:
\begin{align*}
(a)&\quad x_{i+1}\leq x_{i},\\
(b)&\quad x_{i+1}=1+\asc_{i}(x),\mbox{ or}\\
(c)&\quad x_i< x_{i+1}< 1+\asc_{i}(x).
\end{align*}
Now, case $(a)$ occurs exactly when $\addi$ is used on $A^{(i)}$, case $(b)$ when $\addii$ is used, and case $(c)$ when case $\addiii$ is used. 
Furthermore, in all cases, $A^{(i+1)}=\addto{A^{(i)},x_{i+1}}$.

\section{Statistics and distributions}
In this section we show how statistics on the two structures are related under $\bij$.
Many of the definitions concerning ascent sequences were stated in \cite[\S 5]{bcdk}
and we recall them here.

Let $x=(x_1,\ldots , x_n)$ be a sequence of integers.  
For $k\leq n$, define $\asc_k(x)$ to be the number of ascents in the 
subsequence $(x_1, x_2, \ldots, x_k)$. If $x_i<x_{i+1}$, we say that $x_{i+1}$ is an \emph{ascent top}.

Let $\zeros(x)$ be the number of zeros in $x$, and let $\last(x):=x_n$.
A right-to-left maximum of $x$ is an entry $x_i$ that has no larger entry to its right. 
We denote by $\rmax(x)$ the number of right-to-left maxima of $x$.

For sequences $x$ and $y$ of non-negative integers, let $x\oplus y=xy'$, where $y'$ is obtained from $y$
by adding $1+\max(x)$ to each of its letters, and juxtaposition denotes concatenation.
For example $(3,2,0,1,2)\oplus(0,0,1)=(3,2,0,1,2,4,4,5)$.
We say that a sequence $x$ has $k$ {\it{components}} if it is the sum of $k$, but not $k+1$, nonempty nonnegative sequences, 
and write $\comp(x)=k$.

Define $\ascset(x)=\{i: i \in [n-1] \mbox{ and } x_i< x_{i+1}\}$.
We denote by $\hat{x}$ the outcome of the following algorithm;
\medskip\\
\noindent
\begin{minipage}{30em}
\mbox{}{\tt for} $i\in \ascset(x)$: \\
\mbox{}\qquad{\tt for} $j \in [i-1]$: \\
\mbox{}\qquad\qquad{\tt if} $x_j \geq x_{i+1}$ {\tt then} $x_j := x_j+1$
\end{minipage}
\medskip\\
and call $\hat{x}$ the {\it{modified ascent sequence}}.
For example, if $x=(0,1,0,1,3,1,1,2)$ then $\ascset(x)=(1,3,4,7)$ and $\hat{x}=(0,3,0,1,4,1,1,2)$.

Note that the modified ascent sequence $\hat{x}$ has its ascents in the same positions as the original sequence, 
but that the ascent tops in $\hat{x}$ are all distinct.
An ascent sequence $x$ is \emph{self-modified} if $\hat{x}=x$.

Let $\rowsum{i}{A}$ and $\colsum{i}{A}$ be the sum of the elements in row $i$ and column $i$ of $A$, respectively.
Also, let $\flip(A)$ be the reflection of $A$ in its antidiagonal.
Let $\blocks(A)$ be the number of diagonal blocks in the matrix $A$.

\begin{theorem}
\label{thm:rows}
Let $A \in \ourclass{n}$ and $x=\biject(A) \in \Asc{n}$.
Then $$\rowsum{k}{A} = |\{j: \hat{x}_j=k-1\}|.$$ 
\end{theorem}

\begin{proof}
By induction. The result is true for the single matrix $(1) \in \ourclass{1}$.
Let us suppose that the result is true for all matrices $\ourclass{n-1}$ for some $n\geq 2$.
Given $B \in \ourclass{n-1}$, let $x=(x_1,\ldots , x_{n-1})=\biject(B)$
and set $\hat{x}=(\hat{x}_1,\ldots , \hat{x}_{n-1})$.
Let $A = \addto{B,i}$ and $y=(x_1,\ldots , x_{n-1},i) = \biject(A)$.
Furthermore set $\hat{y}=(\hat{y}_1,\ldots , \hat{y}_{n})$.

If $i+1 \leq \ourindex{B}$ then \addi\ applies. 
Thus $\rowsum{k}{A}=\rowsum{k}{B}$ for all $k\neq i+1$, and $\rowsum{i+1}{A}=1+\rowsum{i+1}{B}$. 
Since $i\leq x_{n-1}$ we have that $n-1 \not\in \ascset(x)$, therefore $\hat{y}_j=\hat{x}_j$ for all $j\leq n-1$ and $\hat{y}_n = i$.

If $i=\ourdim{B}$ then \addii\ applies.
In this case $\rowsum{k}{A}=\rowsum{k}{B}$ for all $k\leq \ourdim{B}$ and $\rowsum{\ourdim{B}+1}{A}=1$.
Since $i=1+\asc(x_1,\ldots , x_{n-1})$ we find that
$\hat{y}_j=\hat{x}_j$ for all $j\leq n-1$ and $\hat{y}_n = \max(\hat{x}) +1$.

The remaining case for $\ourindex{B} \leq i< \ourdim{B}$ is dealt with in a similar manner.
\end{proof}

Given a square matrix $A$ and a sequence $x$, define the power series
$$\begin{array}{rclcrcl}
\chi(x,q) &:=& \ds\sum_{i=1}^{|x|} q^{x_i}, & &\overline{\chi}(x,q) &:=& \ds\sum_{x_i \; \text{rl-max}} q^{x_i},\\[2em]
\lambda(A,q)&:=&\ds\sum_{i=1}^{\ourdim{A}} q^{\rowsum{i}{A}}, &&\overline{\lambda}(A,q) &:=& \ds\sum_{i=1}^{\ourdim{A}} A_{i,\ourdim{A}} q^{i-1}.
\end{array}$$

\begin{theorem}\label{thm:stats}
Suppose $A$ is the matrix corresponding to the ascent sequence $x$. Then 
\begin{enumerate}
\item $\zeros(x) = \rowsum{1}{A};$
\item $\last(x)=\ourindex{A}-1;$
\item $\asc(x)=\ourdim{A}-1;$
\item $\rmax(\hat{x})=\colsum{\ourdim{A}}{A};$
\item $\comp(\hat{x})=\blocks(A);$
\item $\chi(\hat{x},q) = \lambda(A,q);$
\item $\overline{\chi}(\hat{x},q) = \overline{\lambda}(A,q)$.
\end{enumerate}
\end{theorem}
\begin{proof}
Most of the results follow from the sequence of rules applied to construct the matrix $A$ from the ascent sequence $x$.

(i)
An entry $x_j=0$ if and only if the corresponding entry of the modified ascent sequence $\hat{x}_j=0$. 
This result now follows from Theorem \ref{thm:rows} with $i=1$.\\

(ii) and (iii) follow directly from Theorem \ref{thm:biject}.\\

(iv) is an immediate consequence of the proof of (vii) below with $q=1$.\\

(v) We now show that $\comp(\hat{x})=\blocks(A)$.
It suffices to prove that $\hat{x} = \hat{y} \oplus \hat{z}$ with $|y|=\ell$ and $|z|=m$ 
iff
$A=\twomatrix A_y & 0 \\ 0 & A_z \ematrix$ with $A_y \in \ourclass{\ell}$ and $A_z \in \ourclass{m}$,
where $\biject(A_y)=y$ and $\biject(A_z)=z$.

Let us assume that $\hat{x}=\hat{y}\oplus \hat{z}$.
The first $\ell$ steps of the construction of $A$ give $A_y$ where $\ourdim{A_y}=\asc(y)+1$.
Next, since $\hat{x}_{\ell +1}=1+\max\{\hat{x}_j : j \leq \ell\}$, the addition rule \addii\ is used,
and we have
$$A'=\twomatrix A_y &0 \\ 0 & 1 \ematrix$$
where the new 1 is in position $(\asc(y)+2,\asc(y)+2)$.
All subsequent additions, $x_j$ for $\ell +1 <j \leq \ell+m$ are such that $\hat{x}_j\geq 1+\asc(y)$,
and so do not affect the first $\asc(y)+1$ rows or columns of $A'$.
Further to this, the construction that takes place for steps $\ell+1 , \ldots ,\ell +m$ 
has the same relative order as the construction of $A_z$. 
This gives 
$$A=\twomatrix A_y & 0 \\ 0 & A_z \ematrix.$$

Conversely assume that $A=\twomatrix B & 0 \\ 0 & C \ematrix$ with $B \in \ourclass{\ell}$ and $C \in \ourclass{m}$ and $n=\ell+m$.
The first $m$ removal operations only affect entries in $C$ since there is at least one non-zero entry in every row and column 
of $C$. 
Thus $\widehat{x}_{\ell +1},\ldots , \widehat{x}_n\geq \ourdim{B}$ and in particular, $\widehat{x}_{\ell+1}=\ourdim{B}$.
Note that the sequence $(x_{\ell+1}-\ourdim{B}, \ldots , x_n - \ourdim{B})=(z_1,\ldots , z_m)$ is an ascent sequence which is $\biject(C)$.
After these removals, we are left with the matrix $B$, and since it is in $\ourclass{\ell}$, 
the values $x_1,\ldots , x_{\ell} <\ourdim{B}$. Let $y_j = x_j$ for all $j\leq \ell$.
Consequently one has $\hat{x}=\hat{y}\oplus \hat{z}$.
\ \\

(vi) is an immediate consequence of Theorem \ref{thm:rows}.\\

Finally, part (vii) is proved by induction as follows.
The result is clearly true for the single matrix $(1) \in \ourclass{1}$.
Assume it is true for all matrices in $\ourclass{n-1}$ for some $n\geq 2$.
Let $B \in \ourclass{n-1}$ with $x'=(x_1,\ldots , x_{n-1}) =\biject(B)$.
Let $A=\addto{B,i}$ with $x=(x_1,\ldots , x_{n})=\biject(A)$.
Then
\begin{eqnarray*}
\overline{\lambda} (A,q) = \left\{
	\begin{array}{ll}
	\overline{\lambda}(B,q) + q^i & \mbox{ if } i \leq \ourindex{B}-1\\[1em]
	q^i+\ds\sum_{j=i+1}^{\ourdim{B}} B_{j,\ourdim{B}} q^{j} & \mbox{ otherwise.}
	\end{array}
	\right.
\end{eqnarray*}
Similarly,
\begin{eqnarray*}
\overline{\chi} (\hat{x},q) = \left\{
	\begin{array}{ll}
	\overline{\chi}(\widehat{x'},q) + q^i & \mbox{ if } i \leq x_{n-1}\\[1em]
	q^i+\ds\sum_{\text{rl-max }\widehat{x'_j}\geq i} q^{\widehat{x'_j}+1} & \mbox{ otherwise.}
	\end{array}
	\right.
\end{eqnarray*}
From the induction hypothesis, for the case $i\leq \ourindex{B}-1=x_{n-1}$, 
we have $\overline{\lambda}(B,q)=\overline{\chi}(\widehat{x'},q)$.
Otherwise, 
$$\sum_{j=i+1}^{\ourdim{B}} B_{j,\ourdim{B}} q^{j}=\sum_{\text{rl-max }\widehat{x'_j}\geq i} q^{\widehat{x'_j}+1}$$ 
since these power series are simply $\overline{\lambda}(B,q)$ and $\overline{\chi}(\widehat{x'},q)$, 
respectively, without the first $i$ powers of $q$.
Thus $\overline{\chi} (\hat{x},q) =\overline{\lambda} (A,q)$.
\end{proof}

The above results, used in conjunction with the $\flip$ operation, allow us to prove the following equidistribution result on ascent sequences.

\begin{theorem}
For all $n\geq 1$,
$\zeros (x)$ and $\rmax (\hat{x})$ are equidistributed on the set $\Asc{n}$.
\end{theorem}

\begin{proof}
\newcommand{\newmap}{\mu}
Given $A \in \ourclass{n}$ with $\rowsum{1}{A}=j$, the matrix $B=\flip (A)$ is such that $\colsum{\ourdim{A}}{B}=j$.
This gives 
$$\sum_{A \in \ourclass{n}} q^{\rowsum{1}{A}} = \sum_{A \in \ourclass{n}} q^{\colsum{\ourdim{A}}{A}}.$$
Using Theorem \ref{thm:stats} (i) and (iv), 
$$\ds\sum_{x \in \Asc{n}} q^{\zeros (x)} = \sum_{x \in \Asc{n}} q^{\rmax (\hat x)},$$
thereby showing that $\zeros (x)$ and $\rmax (\hat{x})$ are equidistributed on $\Asc{n}$.
\end{proof}

In dealing with compositions of an integer, the {\it{number of parts}} in a composition is a natural statistic by which the collection of compositions may be refined. 
The next theorem gives the relation between the  number of non-zero parts in our `matrix composition of the integer $n$' and the ascent sequence to which it corresponds.

\begin{theorem}
Let $x=\bij(A)$ where $A \in \ourclass{n}$. The number of positive entries in $A$ is equal to $n$ less the number of equal adjacent entries in $x$.
\end{theorem}

\begin{proof}
Suppose that $x=(x_1,\ldots , x_n)$. 
Let $A^{(i)}$ be the matrix corresponding to $(x_1,\ldots , x_i)$ and define $N_i$ to be the number of 
positive entries in $A^{(i)}$.
Since $A^{(1)}=(1)$ we have $N_1=1$.
Given $i\geq 2$, if $x_i<x_{i-1}$ then one of the zeros in $A^{(i-1)}$ becomes a one in $A^{(i)}$ so that $N_i=N_{i-1}+1$.
If $x_i=x_{i-1}$ then $\ourvalue{A^{(i-1)}}$ is increased by one to give $A^{(i)}$, so in this case $N_i=N_{i-1}$.
Otherwise $x_i>x_{i-1}$ and a new row and column is inserted into $A^{(i-1)}$ to give $A^{(i)}$, and a 1 is introduced, giving
$N_i=N_{i-1}+1$. These equalities may be summarized by $N_i=N_{i-1}+a_i$ where
$a_i=\indicator{x_i \neq x_{i-1}}$.
So the number of positive entries in $A$ is
$$1+a_2+\ldots +a_n = 1+(n-1)-\sum_i \indicator{x_i=x_{i-1}},$$
which is $n$ less the number of equal adjacent entries in $x$.
\end{proof}

\begin{theorem}
The trace $\text{tr}(A)$ is equal to the number of entries $x_i$ 
in the corresponding sequence $x$ such that
$x_i=\asc_{i}(x)$.
\end{theorem}

\begin{proof}
First note that if $i>1$, then $x_i=\asc_{i}(x)$ if either $x_i=1+\asc_{i-1}(x)$ or if $x_i=x_{i-1}=\cdots=x_{i-j}=1+\asc_{i-j-1}(x)$ for some $j\geq 1$.

Now consider the step-by-step process of building $A$. 
If $x_i=1+\asc_{i-1}(x)$, then the matrix dimension increases, 
and a new entry 1 is inserted at the end of the diagonal. 
If $j>0$, and $x_{i+j}=\cdots=x_{i}=1+\asc_{i-1}(x)$, then the same entry gets increased by one. 

Entries at the diagonal can never decrease, and the two cases above 
are the only times an entry on the diagonal can increase, so the result follows.
\end{proof}

Define a {\it{run}} in the sequence $x=(x_1, \ldots, x_n)$ to be a maximal 
subsequence of adjacent equal elements, that is, 
a subsequence $(x_i, x_{i+1}, \ldots, x_{i+j})$ such that $x_i=x_{i+1}=\cdots =x_{i+j}$, where
$x_{i-1}\not = x_i$ if $i>1$ , and $x_{i+j}\not = x_{i+j+1}$ if $i+j<n$. 
If $x_i=y$, we say the run is a $y$-run.

\begin{theorem}
Let $A\in\ourclass{n}$, and suppose that $x=\bij{(A)}$ is the corresponding sequence. The following three equalities hold.
\begin{enumerate}
\item $A_{1,1}$ equals the length of the starting 0-run.
\item $\ourvalue{A}$ equals the length of the ending $x_n$-run.
\item $A_{\ourdim{A},\ourdim{A}}$ equals the length of the last $y$-run whose first entry $x_i=y$ satisfies $x_i=1+\asc_{i-1}(x)$.
\end{enumerate}
Furthermore, the distribution of all three statistics on matrices are the same, 
as is the distribution of all three statistics on ascent sequences.
\end{theorem}

\begin{proof}
Using the standard method of building the matrix according to the ascent sequence it is straightforward to check that the three equalities hold.

To show that the first two statistics on ascent sequences are equidistributed, a simple bijection can be used. 
Assume that $x$ is of the form $(0^a,y,i^b)$, where the subsequence $y$ starts with 1, and does not end with $i$. 
Map $x$ to $\tilde{x}=(0^b,y,i^a)$. It is obvious that this is a bijection (and also an involution), and that the result follows.

The third statistic also have the same distribution by symmetry --- it is equal to $\flip(A)_{1,1}$.
\end{proof}

\begin{remark}
The observant reader may have noticed that there is a fourth pair missing from the above theorem: the last positive entry in the first row of the matrix, and its counterpart for sequences. The counterpart is the length of a subsequence of zeros, but the rule for deciding which is quite complicated. 
\end{remark}

\begin{conjecture}
For ascent sequences $x$, the distribution of $\zeros(x)$, or equivalently, the distribution of $\rmax(\hat x)$, is the same as the distribution of the length of the first strictly increasing subsequence of $x$.
\end{conjecture}

\section{Binary, positive diagonal, and bidiagonal matrices}
We now turn to some natural subclasses of matrices. 
These are binary matrices, matrices that have no zeros on their diagonal, and bidiagonal matrices.

First, let us note that it is easy to see that the collection of diagonal 
matrices in $\ourclass{n}$ correspond to compositions of the integer $n$. 
Given such a matrix $A=\mbox{diag}(a_0,\ldots , a_k) \in \ourclass{n}$, the corresponding ascent sequence is
$$\biject(A) = (0^{a_0},1^{a_1}, \ldots , k^{a_k}), \mbox{ with } a_0+a_1+\cdots+a_k=n.$$

It is known that the binary matrices in $\ourclass{n}$ correspond to interval orders with no repeated holdings \cite{F85}. 
These are a subclass of interval orders, which were shown in \cite{bcdk} to be in bijection with ascent sequences. 

\begin{theorem}
A matrix $A \in \ourclass{n}$ is binary if and only if the corresponding ascent 
sequence $x=\biject(A)$ contains no two equal consecutive entries.
\end{theorem}

\begin{proof}
Suppose that $A\in \ourclass{n-1}$ is a binary matrix with $x=(x_1,\ldots , x_{n-1})=\bij(A)$.
Let $B=\addto{A,m}$.
Rules \addii\ and \addiii\ introduce a new 1 to a matrix by increasing the dimensions of the matrix by one (and maybe permuting some of the entries).
Rule \addi\ increases the entry of a position in the matrix by $1$. In the case that $0\leq m < \ourindex{A}-1$ the value that is increased by 1 will be 0.
However, if $m=\ourindex{A}-1$ then $\ourvalue{A}>0$ will be increased by one to yield $\ourvalue{B}>1$.
In terms of the corresponding ascent sequence this is equivalent to $x_{n-1}=m$.
\end{proof}

The previous result can be generalized slightly. The proof is similar, and omitted.
\begin{proposition}
Let $A\in\ourclass{n}$ be the matrix corresponding to the ascent sequence $x$. Then the sum 
$\sum_{i,j} \max\{0,(A_{i,j}-1)\}$ equals the number of pairs $(x_i,x_{i+1})$ in $x$ such that $x_i=x_{i+1}$. 
\end{proposition}

Next we classify those matrices in $\ourclass{n}$ that have only positive diagonal entries.
Let us point out that the following class of ascent sequences correspond to permutations that
avoid the pattern $3\overline{1}52\overline{4}$, see \cite[Prop. 9]{bcdk}.

\begin{theorem}
The matrix $A=\invbij(x)$ has only positive entries on the diagonal exactly when the sequence $x$ is self-modified, that is when $x=\hat{x}$.
\end{theorem}

\begin{proof}
Consider the sequence of addition rules used to build $A$. 
If $A$ has no zeros on the diagonal, it means that $\addiii$  was never used. 

This means that for the sequence $x$, for all $i$, $x_{i-1}\geq x_i$ or $x_i=1+\asc_{i-1}(x)$. 
In other words all ascents are maximal. 
This is exactly the condition for a sequence to be self-modified: 
a sequence is not self-modified if and only if there exist $i$ and $j<i$ such that $x_j\geq x_{i+1}$ and $x_i<x_{i+1}$.
\end{proof}

\subsection{Bidiagonal matrices and order-consecutive set partitions}

Consider the subclass $\bidiagclass{n}\subseteq\ourclass{n}$ 
of matrices defined to be the bidiagonal matrices in $\ourclass{n}$. 
It turns out that there is a natural bijection between $k\times k$ matrices in 
$\bidiagclass{n}$ and so called {\it{order-consecutive set partitions}}, \cite{HM1995}, of $[n]$ into $k$ parts. 
A set partition is order-consecutive if the parts $P_1,P_2, \ldots, P_k$ can be ordered as 
$$P_{\pi_{1}},P_{\pi_{2}}, \ldots, P_{\pi_{k}}$$
such that each set $\bigcup_{i=1}^j P_{\pi_i}$ is an interval in $[n]$. For example, 
$$\{\{1,2,3\},\{4,9\},\{5\},\{6,7\},\{8\}\}$$ 
is order-consecutive, for we can order the parts as 
$$\{5\},\{6,7\}, \{8\}, \{4,9\}, \{1,2,3\}.$$ 
The set partition $\{\{1,3\},\{2,4\}\}$ however is not order-consecutive.

An order-consecutive set partition of $[n]$ into $k$ parts can be 
represented as the sequence $1$ to $n$, with $k$ pairs of parenthesis inserted
(see \cite{HM1995}).
For example, $\{\{1,2,3\},\{4,9\},\{5\},\{6,7\},\{8\}\}$ is represented as $(123)(4(5)(67)(8)9)$. 
Note that each pair of parenthesis are placed as close together as possible. 
Thus, $(1(2))$ is not a valid representation --- the proper one for this partition is $(1)(2)$. 
These representations for order-consecutive partitions obey an additional constraint \cite[Lemma 5]{HM1995}:

{\it{Constraint $\ast$:}} If all $)($-pairs are deleted, the remaining pairs are completely nested, 
i.e. removing the numbers we are left with $((\cdots())\cdots)$. 

Given an order-consecutive set partition $P=(P_1,\ldots, P_k)$ of $[n]$, 
let us write $\alpha(P)$ for this representation involving parentheses.
We form a bidiagonal matrix $B=B(P)$ as follows. 
Let $B(P)$ initially be the $k\times k$ matrix with all elements zero except a one at the top left corner. 
Read the sequence $\alpha(P)$ from left to right, starting with 1. 
When reading the sequence, if the next symbol is a number, increase the element in the current position of $B(P)$ by one. 
If it is a parenthesis increase either the row index or column index of $B(P)$ by one, 
whichever allows us to stay on the diagonal and bidiagonal.

For example, the partition above with $\alpha(P)=(123)(4(5)(67)(8)9)$ gives the matrix
$$B(P)=\fivematrix 3&0&0&0&0\\ 0&1&1&0&0\\ 0&0&0&2&0\\ 0&0&0&0&1\\ 0&0&0&0&1\ematrix.$$

\begin{theorem}
There is a one-to-one correspondence between $k\times k$ matrices in $\bidiagclass{n}$ and order-consecutive set partitions of $[n]$ into
$k$ parts. 
\end{theorem}

\begin{proof}
It is clear from above construction that if $P$ is an order-consecutive set partition, then $B(P) \in \bidiagclass{n}$. 
We show it is one-to-one by defining the inverse.
The numbers 1 to $n$ are to be written down in order, with parenthesis interspersed. 
Start by writing $($. 
Next visit the elements in the matrix in order $(1,1), (1,2), (2,2), (2,3), \ldots$. 
If the number $m$ is encountered, write down the next $m$ numbers and then a $|$. End with a $)$. 

Now change each $||$ into $)($. 
Note that there can be no more than two consecutive $|$'s.

To finish, we need to change each remaining $|$ into either $)$ or $($. 
However, using constraint $\ast$ , there is a unique way of doing this.
\end{proof}

\begin{corollary}[\cite{HM1995}]
The number of $k\times k$ bidiagonal matrices in $\ourclass{n}$ is 
$$\sum_{j=0}^{k-1}\binom{n-1}{2k-j-2}\binom{2k-j-2}{j}.$$
Furthermore, from the construction above one may notice that the term 
in the sum counts the number of matrices with exactly $j$ zeros in the diagonal and bidiagonal.
\end{corollary}

\begin{theorem}
The set of ascent sequences $x$ such that $x=\bij (A)$ for $A \in \bidiagclass{n}$ are those sequences $x=(x_1,\ldots , x_n)$
which satisfy
\begin{equation}
x_i\geq \asc_{i}(x)-1,\label{equi1}
\end{equation}
for $1\leq i\leq n$. 
\end{theorem}

\begin{proof}
Induction on $n$. 
The $n=1$ case is trivial, so assume that $A\in\ourclass{n}$ is bidiagonal, 
and that $x=(x_1,x_2, \ldots, x_{n})=\bij(A)$ obeys \eqref{equi1}.

Consider two cases. First assume the last column of $A$ ends with $(0,a)$ for some $a\geq 1$. 
This means that $x$ ends with $x_{n+1-a}=x_{n+2-a}=\cdots =x_{n}=\asc_{n+1-a}(x)=\asc_{n}(x)$. 

Let $y=(x_1, \ldots, x_n, x_{n+1})$ and $B=\invbij(y)$. Consider the three subcases
$x_n<x_{n+1}$, $x_n=x_{n+1}$ and $x_n>x_{n+1}$.

If $x_{n+1}=x_{n}+1$ then $B=\invbij(y)$ is $\twomatrix A&0\\0&1\ematrix$, 
and bidiagonal by the induction hypothesis. 
Also, $x_{n+1}=\asc_{n}(y)+1=\asc_{n+1}(y)$, so $x_{n+1}\geq \asc_{n+1}(y)-1$.

If $x_{n+1}=x_{n}$ then $B$ is $A$ with the entry at position $(\ourdim{A},\ourdim{A})$ 
increased by one, and again bidiagonal. Furthermore, $x_{n+1}=\asc_{n+1}(y)\geq\asc_{n+1}(y)-1$.

If $x_{n+1}=x_{n}-m$ for $m>0$, then $B$ is $A$ with the 0 at position $(\ourdim{A}-m,\ourdim{A})$ 
increased to a 1, and $A$ is bidiagonal if and only if $m=1$. 
Also, $x_{n+1}=\asc_n(x)-m$, so  $x_{n+1}\geq\asc_{n+1}(y)-1$ only for $m=1$.

This proves the theorem in first case. 
The second case, when the last column of $A'$ ends with $(a>0,b>0)$ is handled in a similar way. 
\end{proof}

\section{Some challenging questions}\label{sec:Q}
We end this paper with two challenging questions.

\begin{question}
If $x=\biject (A)$ for some $A \in \ourclass{n}$, then what is the sequence $y=y(x)$ for which $y=\biject (\flip(A))$?
\end{question}

In terms of \tpt -free posets, this question is equivalent to asking for the ascent sequence $y$ that corresponds 
to the dual poset $P^{\star}$, where the poset $P$ is generated by the ascent sequence $x$. 

Adding two upper triangular matrices of the same dimension yields another upper triangular matrix of the same dimension.

\begin{question}
Adding two matrices of the same size is a commutative mapping $\ourclass{n}\times \ourclass{m} \mapsto \ourclass{n+m}$.
How does this operation act on the corresponding ascent sequences?
Furthermore, how does this addition operation act on the corresponding posets?
\end{question}

\end{document}